\documentclass{amsart}

\usepackage[dvipsnames,svgnames,x11names,hyperref]{xcolor}
\usepackage{textgreek,bbm,url,graphicx,verbatim,amssymb,enumerate,stmaryrd,booktabs,lmodern,mathtools,microtype,mathabx}
\usepackage[pagebackref,colorlinks,citecolor=Mahogany,linkcolor=Mahogany,urlcolor=Mahogany,filecolor=Mahogany]{hyperref}
\usepackage[capitalize]{cleveref}
\usepackage[mathscr]{euscript}
\usepackage{tikz,tikz-cd}
\usepackage{pgfplots}

\usepackage{datetime}

\newdateformat{monthdayyeardate}{%
  \monthname[\THEMONTH]~\THEDAY, \THEYEAR}

\usepackage{float}

\usetikzlibrary{matrix,calc,positioning,arrows,decorations.pathreplacing,patterns,arrows}

\newcommand{\N}{\mathbb{N}}

\newcommand{\R}{\mathbb{R}}

\DeclareMathOperator{\Comp}{Comp}

\DeclareMathOperator{\closure}{closure}
\DeclareMathOperator{\CW}{CW}

\newcommand*\interior[1]{#1^{\mathsf{o}}}

\usepackage[pagebackref,colorlinks,citecolor=Mahogany,linkcolor=Mahogany,urlcolor=Mahogany,filecolor=Mahogany]{hyperref}
\usepackage[dvipsnames,svgnames,x11names,hyperref]{xcolor}
\usepackage{url,graphics,verbatim,amssymb,enumerate,stmaryrd,mathtools}
\usepackage{setspace}
\usepackage{microtype}
\usepackage[all,cmtip]{xy}
\usepackage[hmargin=2.75cm]{geometry}

\newtheorem{theorem}{Theorem}[section]
\newtheorem{lemma}[theorem]{Lemma}
\newtheorem{proposition}[theorem]{Proposition}

\theoremstyle{definition}

\newtheorem{remark}[theorem]{Remark}

\title{Complements Of Locally Flat Submanifolds Are Finite CW Complexes}
\author{Andrew Ho}
\date{\monthdayyeardate\today}

\begin{document}

\begin{abstract} We show that if $Y$ is a compact topological manifold and $X$ is a locally flat submanifold, then the complement $Y - X$ is homotopy equivalent to a finite CW complex. This is a direct proof, and does not rely on much of the theory of topological manifolds.\end{abstract}

\maketitle

\section{Introduction}

If $Y$ is a smooth manifold while $X$ is a compact smooth submanifold, then $Y - X$ is homotopy equivalent to a finite CW complex: equipping $Y$ with a Riemannian metric $g$, the normal vector bundle of $X$ yields a closed disk bundle $D$ of $X$ such that $Y - \interior D$ is a compact smooth manifold (thus homotopy equivalent to a finite CW complex, by Theorem 3.5 of \cite{milnormorse}) while $Y - X$ deformation retracts onto $Y - \interior D$. \\
\newline
Attempting to generalize this fact to the case where $X, Y$ are topological manifolds, one requires that $X$ be a \textit{locally flat} (or \textit{tame}) submanifold of $Y$. The reason is that if one takes $X$ to be the Alexander Horned Sphere along with its interior and $Y = S^{3}$, then $\pi_{1}(Y - X)$ is not finitely generated (by \cite{blankfox}), and so $Y - X$ cannot be homotopy equivalent to a finite CW complex. The main result to be proven in this paper is:
\begin{theorem} \label{thm:maintheorem}
Suppose $Y^{n}$ is a compact topological manifold while $X^{d}$ is a compact locally flat submanifold. Then $Y - X$ is homotopy equivalent to a finite CW complex.
\end{theorem}
It is indeed true that if $X$ has a closed disk bundle, then $Y - X$ would be homotopy equivalent to a finite CW complex. This fact will be proven in \cref{sec:assembly}, and is listed below:
\begin{proposition} \label{prop:diskbundle}
Assume the same hypotheses as \cref{thm:maintheorem}. If $X$ has a closed disk bundle in $Y$, then $Y - X$ is homotopy equivalent to a finite CW complex.
\end{proposition}
The proof of \cref{thm:maintheorem} thus seems to reduce to proving \cref{prop:diskbundle} and showing that $X$ has a closed disk bundle. However, closed disk bundles need not exist in general, seen by the following examples:
\begin{itemize}
\item \cite{rourkesanderson} shows that if one takes $Y = S^{29}$, then there exists a compact locally flat submanifold $X$ that does not admit a normal bundle.

\item Theorems 1, 3 of \cite{williambrowder} gives an example where $X$ is locally flat and has a normal bundle, but does not have a closed disk bundle.

\item Theorem 4 of \cite{hirsch} gives examples of locally flat $X^{4} \subset S^{7}$ and $S^{4} \subset M^{7}$ that do not admit closed disk bundles.
\end{itemize}
Even though $X$ need not have a closed disk bundle, the stable uniqueness theorem of \cite{kirbysiebenmann} (also see (B) of \cite{hirschnormal} and Section 5 of \cite{milnornormal}) shows that for some $N \in \N$, $X \times 0$ has a normal bundle in $Y \times \R^{N}$; by this the corollary of Proposition 2.2 of \cite{williambrowder} shows that, after incrementing $N$ by $1$, $X \times 0$ has a closed disk bundle in $Y \times \R^{N}$ where $Y \times \R^{N}$ deformation retracts onto $Y \times [0,1]^{N}$ while fixing $X - 0$. The key of the proof of \cref{thm:maintheorem} thus lies in arguing that the case for the pair $(Y,X)$ reduces to that for the pair $(Y \times [0,1]^{N}, X \times 0)$. This is done by using the sum formula of Wall's finiteness obstruction discussed in \cite{siebenmannthesis}, where regular open neighbourhoods discussed in \cite{siebenmann} are used to show that $Y - X$ is finitely dominated, allowing the consideration of Wall's finiteness obstruction. 

\begin{remark}\label{rem:altproof} Our proof of \cref{thm:maintheorem} avoids a potential proof relying on several cases and much of the theory of topological manifolds. An outline of this proof is:
\begin{itemize}
\item[(Case 1)] When $n - d \leq 2$: we use the same argument as presented in \cref{sec:codimleq2}.

\item[(Case 2)] When $n - d \geq 3$ and $n \geq 5$: we show that $X$ admits a mapping cylinder neighborhood $\text{Cyl}(f) = ((W \times [0,1]) \sqcup X)/\sim$ with $(w,1) \sim f(x)$ where $f: W \to X$. For $n \geq 6$, this is already done in \cite{pedersen} and \cite{endone}, and that of \cite{endone} can supposedly be generalized to $n \geq 5$ using \cite{endthree}.

\item[(Case 3)] When $n = 4$ and $d = 1$: We use \cite{freedmanquinn} to obtain a normal vector bundle of $X$ in $Y$, by which we obtain a closed disk bundle of $X$ in $Y$ and apply \cref{prop:diskbundle}.
\end{itemize}
\end{remark}

Henceforth we assume that $Y$ is connected, as it only has finitely many components by compactness. To be able to use the sum formula in \cite{siebenmannthesis}, we require $n - d \geq 3$ so that the inclusion $Y - X \xhookrightarrow{} Y$ induces an isomorphism of fundamental groups. Thus, in \cref{sec:prelim}, we first show that \cref{thm:maintheorem} holds for $n - d \leq 2$, and then discuss the assumption that $n - d \geq 3$ to prepare for the proof of \cref{thm:maintheorem} under this assumption, which the later sections of the paper concern. It should be noted that the proof of \cref{prop:diskbundle} in \cref{sec:assembly} does not rely on any consideration of codimension, so we will be able to use \cref{prop:diskbundle} in the proof of \cref{thm:maintheorem} in low codimension. \\
\newline
It is worth noting that the aforementioned stable existence theorem requires a lack of boundary of $X, Y$. We will thus first assume until \cref{sec:assembly} that $\partial X = \emptyset$, and in \cref{sec:assembly}, we will address the case where $\partial X \neq \emptyset$ and $\partial X \subseteq \partial Y$.

\subsection*{Acknowledgements}
We thank Alexander Kupers for helpful discussions and comments. 

\tableofcontents

\section{Proof for $n - d \leq 2$ and preparation for $n - d \geq 3$}\label{sec:prelim}

\subsection{Proof of Theorem 1.1 when $1 \leq n - d \leq 2$} \label{sec:codimleq2}

When $n - d = 1$, the fact that $X$ is locally flat in $Y$ implies that $X$ admits a bicollar, by Theorem 3 of \cite{brownnormal}. This bicollar is then a closed disk bundle, since $[-1,1]$ is the $1$-dimensional closed unit disk, and so \cref{prop:diskbundle} finishes the proof of \cref{thm:maintheorem}. \\
\newline
For the case where $n - d = 2$: it suffices to show that $X$ admits a normal vector bundle in $Y$, since that normal vector bundle gives us a closed disk bundle, by which \cref{prop:diskbundle} proves \cref{thm:maintheorem}. If $n \neq 4$, then the existence of a normal vector bundle is given by Theorem A of \cite{kirbysiebenmannnormal}; if $n = 4$, then the existence of a normal vector bundle is given by Section 9.3 of \cite{freedmanquinn} (Freedman-Quinn defines a ``normal bundle" as a normal vector bundle).

\subsection{Proof of Theorem 1.1 when $n - d = 0$}

In this case, $\partial X$ is compact locally flat submanifold of $Y$, with codimension $1$. Since we have proven Theorem 1.1 for codimension $1$, we thus deduce that $Y - \partial X$ is homotopy equivalent to a finite CW complex, so that the path components of $Y - \partial X$ are all homotopy equivalent to finite CW complexes. Since $Y - X$ and $X$ only have finitely many components, let us denote the path components of $Y - X$ as $C_{1},...,C_{N}$ and denote the path components of $X$ as $D_{1},...,D_{M}$. It then suffices to show that $\interior D_{1},..., \interior D_{M}, C_{1},...,C_{N}$ are the path components of $Y - \partial X$, in which case we would be able to conclude that $Y - X$ is homotopy equivalent to a finite CW complex. \\
\newline
We first note that the union of $\interior D_{1},..., \interior D_{M}, C_{1},...,C_{N}$ is $Y - \partial X$, and this union is disjoint. Additionally, these sets are open in $Y - \partial X$ due to the identity $n = d$ along with local flatness of $X$. As $\interior D_{1}, ..., \interior D_{M}$ are the path components of the interior of $M$, we deduce that $\interior D_{1},..., \interior D_{M}, C_{1},...,C_{N}$ are path connected in $Y - \partial X$, and so we are left with showing that they are maximally path-connected in $Y - \partial X$. \\
\newline
Since $C_{1},...,C_{N}$ are the path components of $Y - X$ while $\interior D_{1},..., \interior D_{M}$ are the path components of the interior of $X$ while $Y - \partial X$ is the disjoint union of these $M + N$ sets, it suffices to show that given $p \in C_{i}$ and $q \in D_{j}$, there cannot exist a path in $Y - \partial X$ between $p$ and $q$. Indeed, supposing for contradiction that $\gamma$ were such a path with $\gamma(0) = p$ and $\gamma(1) = q$, noting that $Y - X$ and the interior of $X$ are open in $Y$, by connectedness of $[0,1]$ we deduce that $\gamma$ must intersect $\partial X$ somewhere, a contradiction.

\subsection{$n - d \geq 3$: first facts} We show preliminary facts for the proof of \cref{thm:maintheorem} when $n - d \geq 3$. $Y$ and $X$ are given in more generality here, with the compactness assumption removed.

\begin{lemma}\label{lem:codimtwoconnected} Suppose $Y^{n}$ is a connected topological manifold and $X^{d}$ is a locally flat submanifold. If $n - d \geq 2$, then the inclusion $Y - X \xhookrightarrow{} Y$ is $0$-connected. 
\begin{proof}
Since $Y$ is connected, the induced map $\pi_{0}(Y - X) \to \pi_{0}(Y)$ is immediately surjective. Injectivity of this map follows from Theorem 2 of \cite{eilenburgwilder}, where we take $([0,1], \{0,1\}, Y - X, Y)$
in our notation as $(B,Z,A,X)$ in the notation of \cite{eilenburgwilder}; in our notation, we have $\overline{Y - X} = Y$ since $X$ has positive codimension (one applies either invariance of domain or local flatness), and the relevant hypothesis of Theorem 2 of \cite{eilenburgwilder} also follows by local flatness of our $X \subseteq Y$ along with the fact that $n - d \geq 2$. We also want to show that the induced map $\pi_{1}(Y - X) \to \pi_{1}(Y)$ is surjective, and this also follows from Theorem 2 of \cite{eilenburgwilder} with $(S^{1} \times [0,1], S^{1}, Y - X, Y)$ in our notation taken as $(B,Z,A,X)$; the same comment as our first application of Theorem 2 of \cite{eilenburgwilder}.
\end{proof}
\end{lemma}

\begin{lemma}\label{lem:codimthreeconnected} Suppose $Y^{n}$ is a connected topopological manifold and $X^{d}$ is a locally flat submanifold. If $n - d \geq 3$, then the inclusion $Y - X \xhookrightarrow{} Y$ induces an isomorphism of fundamental groups.
\begin{proof}
The inclusion is $0$-connected by \cref{lem:codimtwoconnected}, so it suffices to show that the induced surjection $\pi_{1}(Y - X) \to \pi_{1}(Y)$ is injective. We apply Theorem 2 of \cite{eilenburgwilder}, taking $(D^{2}, S^{1}, Y - X, Y)$
in our notation as $(B,Z,A,X)$ in the notation of \cite{eilenburgwilder}; in our notation we have $\overline{Y - X} = Y$, and the relevant hypothesis follows by local flatness and $n - d \geq 3$.
\end{proof}
\end{lemma}

\section{$Y - X$ is finitely dominated}\label{sec:findom}

To show that $Y - X$ is finitely dominated, we use regular open neighbourhoods developed in \cite{siebenmann}, \cite{sgh73}, and \cite{sgh74}. In particular, we note the following terminology used:
\begin{itemize}
\item \textit{$I$-compressibility} and the \textit{axiom $I$-$\Comp(Y,X)$}, introduced in \cite{siebenmann} right after Theorem 1.3.

\item \textit{A pinch of $A$ in $U$ with respect to $(Y,X)$}, introduced in Definition 2.1 of \cite{siebenmann}.

\item \textit{The pinching axiom $P(x)$}, introduced in Definition 2.1 and Theorem 2.2 of \cite{siebenmann}.

\item \textit{An anti-$I$-nest $\{F_{n}\}$}, appearing in \cite{sgh73} and \cite{sgh74} only. It is introduced in Definition 3.1 of \cite{sgh73} and originally named ``\textit{anti-$I$-gigogne}".
\end{itemize}
To show that $Y - X$ is finitely dominated, we will first show:
\begin{proposition}\label{prop:icomp} The axiom $I$-$\Comp(Y,X)$ holds.
\end{proposition}
Since $Y$ is locally compact and $X$ is compact in $Y$, by Proposition 5.7 of \cite{sgh73}, it suffices to show that the pinching axiom $P(x)$ holds for every $x \in X$; recall the pinching axiom $P(x)$ means that for every neighborhood $x \in U \subseteq Q$, there existence a pinch of $x$ in $U$ with respect to $(Y,X)$. In order to show this, we consider a simplification below.
\begin{lemma} \label{lem:reducetolocalpinch}
To prove \cref{prop:icomp}, it suffices to show existence a pinch of $0$ in $B(0,10)$ with respect to $(\R^{n}, \R^{d} \times 0)$.
\begin{proof}
Fixing $x \in X$ and considering a neighborhood $x \in U \subseteq Y$, to prove \cref{prop:icomp} we want the existence of a pinch of $x$ in $U$ with respect to $(Y,X)$. By local flatness of $X$, we choose an open neighbourhood $x \in V \subseteq Y$ and a homeomorphism $\varphi: V \to \R^{n}$ such that $\varphi(V \cap X) = \R^{d} \times 0$. \\
\newline
A pinch of $x$ in an open subset of $U$ (containing $x$) will also be a pinch of $x$ in $U$. The first implication is that we can take $U$ to be an open neighborhood of $x$. Furthermore, by intersecting $U$ with $V$ and the preimage of an open ball in $\R^{n}$ under $\varphi$, we can assume that $U \subseteq V$ and that $\varphi(U)$ is an open ball in $\R^{n}$. By scaling and translating $\varphi$, we can also assume that $\varphi(U) = B(0,10)$ and that $\varphi(x) = 0$, where $B(0,10)$ denotes the open ball centered at $0$ with radius $10$. \\
\newline
We then see that it suffices to show the existence a pinch of $0$ in $B(0,10)$ with respect to $(\R^{n}, \R^{d} \times 0)$. Given such a pinch $F_{t}: \R^{n} \to \R^{n}$, the homotopy $\varphi^{-1} \circ F_{t} \circ \varphi: V \to V$ would be a pinch of $x$ in $\varphi^{-1}(B(0,10)) = U$ with respect to $(V,X \cap V)$, and since (by the definition of a pinch) this homotopy is supported in the compact subset $\varphi^{-1}(\overline{B(0,20)})$ where $\overline{B(0,20)}$ denotes a closed ball, we can extend it outside of $V$ to obtain a pinch of $X$ in $U$ with respect to $(Y,X)$.
\end{proof}
\end{lemma}
Due to \cref{lem:reducetolocalpinch}, the proof of \cref{prop:icomp} reduces to the proof of \cref{lem:localpinch}.
\begin{proposition}\label{lem:localpinch} There exists a pinch of $0$ in $B(0,10)$ with respect to $(\R^{n}, \R^{d} \times 0)$.
\begin{proof}
We first define the homotopy $H_{t}: [0,\infty) \to [0,\infty)$ by letting $H_{t}$ fix $[1,\infty)$ and having its graph on $[0,1]$ be the linear interpolation of the four points
\begin{align*}
(0,0), (\frac{1-t}{4}, \frac{1-t}{4}), (\frac{1}{2}, \frac{1 - t}{2}), (1,1)
\end{align*}
$H_{0}$ is the identity, and since $\frac{1-t}{4} = 0$ as when $t=1$, the graph of $H_{1}$ on $[0,1]$ is the linear interpolation of $(0,0), (\frac{1}{2},0), (1,1)$. Addtionally, note that $H_{t} \leq H_{0}$. Figure $1$ illustrates the described homotopy.


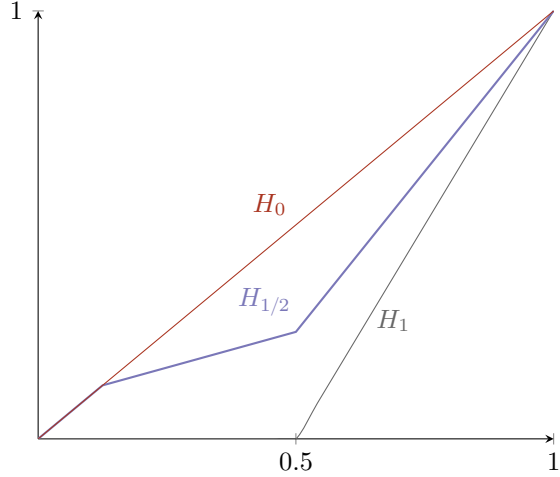
\begin{figure}[H]
\centering
\begin{tikzpicture}[
  declare function={
            func1(\x)= (\x<=0.5) * 0 + (\x > 0.5) * ((2 * \x) - 1);
            funchalf(\x)= (\x <= 1/8) * (x)   +
              and(\x >= 1/8, \x <= 1/2) * (1/3 * \x + 1/12)     +
              (\x >= 1/2) * (3/2 * \x - 1/2);
            func0(\x)= \x;
        }
        ]
        \begin{axis}[
        axis x line=middle, axis y line=middle,
        ymin=0, ymax=1, ytick={0,1},
        xmin=0, xmax=1, xtick = {0,0.5, 1}
        ]
        \addplot[black!60!white, domain=0:1, smooth]{func1(x)} 
        node [pos=0.5, right] {$H_{1}$};
        \addplot[Periwinkle, domain=0:1, thick]{funchalf(x)}
         node [pos=0.4, above left] {$H_{1/2}$};
        \addplot[Mahogany, domain=0:1, smooth]{func0(x)}
        node [pos=0.5, above left] {$H_{0}$};
\end{axis}
\end{tikzpicture}
\caption{The homotopy $H_{t}$.}
	\label{fig:rat}
\end{figure}

Let $\psi: \R^{d} \to [0,1]$ be a continuous function that is constant $1$ in the closed ball $\overline{B(0,1)}$ and supported in the closed ball $\overline{B(0,2)}$. We now define the homotopy $F_{t}: \R^{n} \to \R^{n}$ with
\begin{align*}
F_{t}(u,v) = \begin{cases} (u, \frac{H_{t \psi(u)}(\lvert v \rvert)}{\lvert v \rvert} v) & v \neq 0\\
(u,0) & v = 0 \end{cases}
\end{align*}
For $(u,v) \in \R^{n}$ where $u \in \R^{d}$ and $v \in \R^{n - d}$; continuity is seen by 
noting that $H_{t \psi(u)}(\lvert v \rvert) \leq H_{0}(\lvert v \rvert) = \lvert v \rvert$. \\
\newline
$F_{0}$ is the identity. $F_{1}$ maps $B(0,1) \times B(0,\frac{1}{2})$ into $\R^{d} \times 0$ since $\psi$ is constant $1$ on $\overline{B(0,1)}$ while $H_{1}$ is constant $0$ on $[0,\frac{1}{2}]$. We also see that $F_{t}(u,v) = (u,v)$ for all $t \in [0,1]$ and $(u,v) \in \R^{n} - B(0,10)$, since we would either have $\lvert u \rvert > 2$ so that $\psi(u) = 0$ hence $H_{t\psi(u)}(\lvert v \rvert) = H_{0}(\lvert v \rvert) = \lvert v \rvert$, or $\lvert v \rvert > 1$ so that $H_{t\psi(u)}(\lvert v \rvert) = \lvert v \rvert$ since $H_{t\psi(u)}$ fixes $[1,\infty)$. Finally we note that for each $\tau \in [0,1)$, the family $\{F_{t}: t \in [0,\tau]\}$ is an isotopy that fixes $\R^{d} \times B(0,\frac{1-\tau}{4})$: for $t \in [0,\tau]$ and $\lvert v \rvert < \frac{1-\tau}{4}$ we have $\frac{1-\tau}{4} \leq \frac{1-t}{4} \leq \frac{1-t \psi(u)}{4}$ thus $H_{t\psi(u)}(\lvert v \rvert) = \lvert v \rvert$, and the family $\{H_{s}: s \in [0,\tau]\}$ is an isotopy. This completes the proof.
\end{proof}
\end{proposition}

\noindent We have thus proven \cref{prop:icomp}, showing that the axiom $I$-$\Comp(Y,X)$ holds. By the definition of $I$-$\Comp(Y,X)$, we choose neighborhoods $E_{0}, E_{1}$ of $X$ such that $E_{0} \subseteq E_{1}$ and that $E_{0}$ is $I$-compressible to $X$ in $E_{1}$. Then since $I$-$\Comp(Y,X)$ holds, by Proposition 3.4 of \cite{sgh73} we choose an anti-$I$-nest $\{F_{n}\}$ such that, denoting $B = \cap_{n=0}^{\infty} F_{n}$, we have
\begin{align*}
E_{0} \subseteq B \subseteq F_{0} \subseteq E_{1}.
\end{align*}
To prepare for the proof that $Y - X$ is finitely dominated, we first show the below fact.
\begin{lemma}\label{lem:precompact}
$Y - B$ is precompact in $Y - X$. In other words, the relative closure $\closure_{Y - X}(Y - B)$ is compact in $Y - X$. 
\begin{proof}
It suffices to show that $\overline{Y - B}$ (the closure in $X$) is contained in $Y - X$. The reason is that we would have $\closure_{Y - X}(Y - B) = \overline{Y - B} \cap (Y - X) = \overline{Y - B}$, and since $\overline{Y - B}$ is compact in $Y$ (as $Y$ is compact), we deduce that $\overline{Y - B}$ is compact in $Y - X$, which is what we wanted to conclude. Indeed,
\begin{align*}
\overline{Y - B} &= Y - \interior B &&\text{(The interior is with respect to $Y$)} \\
&\subseteq Y - X &&\text{(Since $E_{0} \subseteq B$ while $E_{0}$ is a neighborhood of $X$)}
\end{align*}
As required. \end{proof}
\end{lemma}

Since $Y$ is metrizable, $X$ as a closed subset is also a $G_{\delta}$ set. So by Example 3.3.(ii) and Theorem 3.6 of \cite{sgh73}, the inclusion $Y - B \xhookrightarrow{\iota} Y - X$ is a homotopy equivalence. With this result, we finally show that $Y - X$ is finitely dominated.

\begin{proposition}\label{prop:findom}$Y - X$ is finitely dominated.
\begin{proof}
$Y - X$ is a topological manifold, and thus an absolute neighborhood retract (ANR); see \cite{milnorhomtype} (Corollary 1) and \cite{hanner} (Theorem 3.3). Due to \cref{lem:precompact}, we can apply Lemma 1.2 of \cite{sgh74} to choose a finite CW complex $L$ and continuous maps $i: Y - B \to L$ and $p: L \to Y - X$ where the following diagram commutes up to homotopy:

\[
\begin{tikzcd}
    Y - B \arrow[hookrightarrow]{rr}{\iota} \arrow[swap]{dr}{i} & & Y - X \\[10pt]
    & L \arrow{ur}{p}
\end{tikzcd}
\]

\noindent Letting $F: Y - X \to Y - B$ denote a homotopy inverse of $\iota$, we then have
\begin{align*}
p \circ (i \circ F) = (p \circ i) \circ F \cong \iota \circ F \cong 1_{Y-X}.
\end{align*}
Since $L$ is a finite CW complex, we conclude that $Y - X$ is finitely dominated.
\end{proof}
\end{proposition}

\section{Reduction to $(Y \times \R^{N}) - (X \times 0)$} \label{sec:reduction}
In this section we aim to show the following:

\begin{proposition}\label{prop:reductionmain} For every $N \in \N$, if $(Y \times \R^{N}) - (X \times 0)$ is homotopy equivalent to a finite CW complex then $Y - X$ is homotopy equivalent to a finite CW complex.
\end{proposition}

\cref{prop:reductionmain} follows by \cref{lem:reductiononedimsimplify} and \cref{prop:reductiononedim}:

\begin{lemma}\label{lem:reductiononedimsimplify}
If $(Y \times \R) - (X \times 0)$ being homotopy equivalent a finite CW complex implies $Y - X$ being homotopy equivalent a finite CW complex, then \cref{prop:reductionmain} holds.

\begin{proof}
We will use the assumption to prove \cref{prop:reductionmain} via induction on $N \in \N$. The base case is precisely the assumption, so we now verify the induction step. Supposing that $(Y \times \R^{N+1}) - (X \times 0)$ is homotopy equivalent to a finite CW complex $K$, we define the compact topological manifolds $Y',X'$ with $X'$ being a locally flat submanifold $Y'$, by the following:
\begin{align*}
Y' = Y \times [-1,1]^{N} \hspace{0.5cm} X' = X \times 0.
\end{align*}
$Y' \times \R - (X' \times 0)$ is homotopy equivalent to $(Y \times \R^{N+1}) - (X \times 0)$ by considering a deformation retract of $\R^{N}$ onto $[-1,1]^{N}$. Therefore, $Y' \times \R - (X' \times 0)$ is homotopy equivalent to $K$. \\
\newline
Since the dimension of $Y'$ is $n + N$ and the dimension of $X'$ is $d$, the codimension of $X'$ in $Y'$ is $n + N - d \geq n - d = 3$. So since $Y' \times \R - (X' \times 0)$ is homotopy equivalent to $K$, we apply the base case to deduce that $Y' - X'$ is homotopy equivalent to a finite CW complex. Now again considering a deformation retract of $\R^{N}$ onto $[-1,1]^{N}$, we see that $Y' - X'$ is homotopy equivalent to $(Y \times \R^{N}) - (X \times 0)$. This means that $(Y \times \R^{N}) - (X \times 0)$ is homotopy equivalent to a finite CW complex, and so applying the induction hypothesis, we conclude that $Y - X$ is homotopy equivalent to a finite CW complex.
\end{proof}

\end{lemma}

\begin{proposition}\label{prop:reductiononedim}
If $(Y \times \R) - (X \times 0)$ is homotopy equivalent to a finite CW complex, then $Y - X$ is homotopy equivalent to a finite CW complex.
\begin{proof}
By the hypothesis, $(Y \times [-1,1]) - (X \times 0)$ is homotopy equivalent to a finite CW complex. Define 
\begin{align*}
U = (Y \times (0,1]) \cup ((Y - X) \times 0) \hspace{0.5cm} V = (Y \times [-1,0)) \cup ((Y - X) \times 0)
\end{align*}
Where $(Y - X) \times 0$ is a topological manifold with dimension $n$. Note that $U,V$ are homotopy equivalent to $Y$, as they both deformation retract onto $Y \times 1$. Also note that $U \cup V = (Y \times [-1,1]) - (X \times 0)$, which is a topological manifold, and connected by \cref{lem:codimtwoconnected}. By the hypothesis, $U \cup V$ is homotopy equivalent to a finite CW complex. We also note that
$U \cap V = (Y - X) \times 0$, so that the inclusion $U \cap V \xhookrightarrow{} Y \times 0$ induces an isomorphism of fundamental groups by \cref{lem:codimthreeconnected}. Finally, since $Y \times 0 \xhookrightarrow{} Y \times [-1,1]$ induces an isomorphism of fundamental groups while the inclusion $U \cup V \xhookrightarrow{} Y \times [-1,1]$ induces an isomorphism of fundamental groups by \cref{lem:codimthreeconnected}, by considering the following commutative diagram of inclusions
\[
\begin{tikzcd}
    U \cap V \arrow[hookrightarrow]{rr}{} \arrow[hookrightarrow]{dr}{} & & Y \times [-1,1]  \\[10pt]
    & U \cup V \arrow[hookrightarrow]{ur}{}
\end{tikzcd}
\]
We deduce that the inclusion $U \cap V \xhookrightarrow{} U \cup V$ induces an isomorphism of fundamental groups.  \\
\newline
By Theorems 2.4-2.5 of \cite{functorial}, we choose CW approximations $\CW(U \cup V), \CW(U), \CW(V), \CW(U \cap V)$ such that each of the four topological spaces is weakly homotopy equivalent to their CW approximation, that the other three CW approximations are subcomplexes of $\CW(U \cup V)$, that $\CW(U \cap V) = \CW(U) \cap \CW(V)$, and that we obtain the following commutative diagram:
\[
\begin{tikzcd}
    U \cap V \arrow{r}{} \arrow[hookrightarrow]{d}{}& \CW(U \cap V) \arrow[hookrightarrow]{d}{} \\[10pt]
    U \cup V  \arrow{r}{} & \CW(U \cup V)
\end{tikzcd}
\]
$U, V, U \cap V, U \cup V$ are connected topological manifolds up to homotopy, while topological manifolds are homotopy equivalent to CW complexes (Corollary 1 of \cite{milnorhomtype}), so by Whitehead's theorem, all weak homotopy equivalences here are homotopy equivalences. In particular, as $U \cap V = (Y - X) \times 0$ is finitely dominated by \cref{prop:findom}, it suffices to show that $\sigma(\CW(U) \cap \CW(V))$ vanishes. \\
\newline
Since $U,V$ are homotopy equivalent to the compact connected topological manifold $Y$, where compact topological manifolds are homotopy equivalent to finite CW complexes by \cite{kirbysiebenmanncompact}, $U,V$ are in particular finitely dominated so that their finiteness obstructions are well defined, with $\sigma(\CW(U)) = \sigma(\CW(V)) = 0$. Additionally, since $U \cup V$ is homotopy equivalent to a finite CW complex, its finiteness obstruction is also well-defined with $\sigma(\CW(U \cup V)) = 0$. Thus by the sum formula for finiteness obstruction (Theorem 6.5 of \cite{siebenmannthesis}), we have
\begin{align*}
0 = 0 + 0 - j_{\ast}\sigma(\CW(U) \cap \CW(V))
\end{align*}
Where $j_{\ast}$ is induced by the inclusion $j: \CW(U) \cap \CW(V) \xhookrightarrow{} \CW(U \cup V)$. But $j$ induces an isomorphism of fundamental groups due to the previous commutative diagram where the inclusion $U \cup V \xhookrightarrow{} U \cap V$ induces an isomorphism of fundamental groups, and so we necessarily have $\sigma(\CW(U) \cap \CW(V)) = 0$.
\end{proof}
\end{proposition}

\section{Assembly of Proof for \cref{thm:maintheorem}} \label{sec:assembly}

\noindent We finish the proof of \cref{thm:maintheorem}, by proving \cref{prop:diskbundle} and considering the case where $\partial X \neq \emptyset$ (and $\partial X \subseteq \partial Y$) and showing that $X$ still admits a normal bundle in $Y \times \R^{N}$ for some $N \in \N$. Recall that \cref{prop:diskbundle} states that the existence of a closed disk bundle $B$ for $X$ implies that $Y - X$ is homotopy equivalent to a finite CW complex.

\subsection{Proof of \cref{prop:diskbundle}} By local triviality of the closed disk bundle $B$, the complement $Y - \interior B$ is a compact manifold, thus homotopy equivalent to a finite CW complex. Therefore, it suffices to show that $Y - \interior B$ and $Y - X$ are homotopy equivalent. \\
\newline
We note that $\partial B$ is a neighbourhood deformation retract (NDR) of $Y - \interior B$: local triviality of $B$ along with Satz 2 of \cite{dold} shows that $(Y - \interior B, \partial B)$ is closed cofibration, by which we apply Proposition 5.4.4 of \cite{tomdieck}. Additionally, attaching $B - X$ to $Y - \interior B$ along $\partial B$ yields $Y - X$, giving us the pushout diagram:

\[
\begin{tikzcd}[font=\large]
    \partial B \arrow[hookrightarrow]{r}{\text{NDR}} \arrow[hookrightarrow]{d}{}& Y - \interior B \arrow[hookrightarrow]{d}{} \\[10pt]
    B - X \arrow[hookrightarrow]{r}{} & Y - X
\end{tikzcd}
\]

\noindent We now show that the inclusion $\partial B \xhookrightarrow{} B - X$ is a homotopy equivalence. Recalling that $B$ is a closed disk bundle over $X$, we let $k$ denote the dimension of this closed disk.

\[
\begin{tikzcd}[font=\large] 
S^{k - 1} \arrow[hookrightarrow]{r}{\simeq} \arrow{d}{}& D^{k} - \{0\} \arrow{d}{} \\[10pt] \partial B \arrow[hookrightarrow]{r}{} \arrow{d}{\pi_{1}} & B - X \arrow{d}{\pi_{2}} \\ [10pt] X \arrow{r}{=} & X
\end{tikzcd}
\]  

\noindent Since $B$ is a $D^{k}$ bundle over $X$, we see that $\partial B$ and $B - X$ are $S^{k-1}$ and $D^{k} - \{0\}$ bundles over $X$ respectively, where the inclusion $S^{k-1} \xhookrightarrow{} D^{k} - \{0\}$ is a homotopy equivalence. By local triviality of $\partial B$ and $B - X$, along with paracompactness of $X$ and the inclusion $S^{k-1} \xhookrightarrow{} D^{k} - \{0\}$ being a homotopy equivalence, we choose a numerable open cover $\{U_{j}: j \in J\}$ of $B$ such that each induced map $\pi_{1}^{-1}(U_{j}) \to \pi_{2}^{-1}(U_{j})$ is a fiberwise homotopy equivalence. By Theorem 13.3.3 of \cite{tomdieck}, we thus deduce that the inclusion $\partial B \xhookrightarrow{} B - X$ is a homotopy equivalence. \\
\newline
Having shown that the inclusion $\partial B \xhookrightarrow{} B - X$ is a homotopy equivalence, we then deduce that $Y - \interior B \xhookrightarrow{} Y - X$ is a homotopy equivalence by applying Proposition 5.3.3 of \cite{tomdieck} with the following commutative diagram in mind:

\[
\begin{tikzcd}[font=\large] 
    \partial B \arrow[hookrightarrow]{rrr}{\text{NDR}} \arrow[hookrightarrow]{ddd}{} &  & & Y - \interior B \arrow[hookrightarrow]{ddd}{} \\[10pt]
    & \partial B \arrow{lu}{=} \arrow[hookrightarrow]{r}{\text{NDR}} \arrow{d}{=} & Y - \interior B \arrow{ru}{=} \arrow{d}{=} & \\[10pt] 
    & \partial B \arrow[hookrightarrow]{ld}{\simeq} \arrow[hookrightarrow]{r}{\text{NDR}} & Y - \interior B \arrow[hookrightarrow]{rd} & \\[10pt] 
    B - X \arrow[hookrightarrow]{rrr}{} &  & & Y - X
\end{tikzcd}
\]

\noindent where the three labelled diagonal maps are homotopy equivalences, and the inner and outer squares are homotopy cocartesian. Having shown that $Y - \interior B \xhookrightarrow{} Y - X$ is a homotopy equivalence, we conclude that $Y - X$ is homotopy equivalent to a finite CW complex, completing the proof of \cref{prop:diskbundle}.

\subsection{Consideration of Boundary}

Here, we suppose $\partial X \neq \emptyset$ and $\partial X \subseteq \partial Y$. Since $Y$ is compact while $X$ is compact and locally flat in $Y$, we see that $\partial X$ is locally flat in $\partial Y$ with $\partial X, \partial Y$ both being closed topological manifolds. Their lack of boundaries allow us to apply the stable existence theorem from \cite{kirbysiebenmann}, so that $\partial X \times 0$ admits a normal microbundle in $\partial Y \times \R^{N}$ for some $N \in \N$. \\
\newline
We apply generalization A.4(b) of \cite{kirbysiebenmann} with $(Y \times \R^{N}, \{\partial Y \times \R^{N}\}, X \times 0)$ in our notation as $(W, \{W_{\alpha}\}, M)$ in the notation of \cite{kirbysiebenmann}, to obtain a map $f: \partial Y \times \R^{N} \times [0,\infty) \to Y \times \R^{N}$ such that $\partial X \times 0 \times [0,\infty) = f^{-1}(X \times 0)$. This thus allows us to extend the normal bundle of $\partial X \times 0$ in $\partial Y \times \R^{N}$ to a normal bundle of $\partial X \times 0$ in $Y \times \R^{N}$; applying the stable existence theorem from \cite{kirbysiebenmann} again, after possibly incrementing $N$, we obtain a normal bundle of $X \times $ in $Y \times \R^{N}$, which is what we wanted. We have thus dealt with the case $\emptyset \neq \partial X \subseteq \partial Y$.


\end{document}